\documentclass[12pt]{article}
\usepackage{epsfig}
\usepackage{amsmath}
\usepackage{amssymb}
\usepackage{eucal}

\newtheorem{thm}{Theorem}[section] 
\newtheorem{lemma}[thm]{Lemma}

\newcommand{\qed}{\hfill\mbox{\raggedright\rule{.07in}{.1in}}
  \vspace{1ex}} 

\newcommand{\Section}[1]{\section{#1} \setcounter{equation}{0}}

\title{Realizability of the normal form for the triple-zero nilpotency in a class of delayed nonlinear oscillators}
 
\author{Victor G. LeBlanc\\Department of
  Mathematics and Statistics\\University of Ottawa\\Ottawa, 
ON K1N 6N5\\CANADA}
\date{\today}

\begin{document}

\maketitle

\begin{abstract}
The effects of delayed feedback terms on nonlinear oscillators has been extensively studied, and have important applications in many areas of science and engineering.  We study a particular class of second-order delay-differential equations near a point of triple-zero nilpotent bifurcation.  Using center manifold and normal form reduction, we show that the three-dimensional nonlinear normal form for the triple-zero bifurcation can be fully realized at any given order for appropriate choices of nonlinearities in the original delay-differential equation.
\end{abstract}

\pagebreak
\Section{Introduction}
Delay-differential equations are used as models in many areas of science, engineering, economics and beyond \cite{BBL, HFEKGG,Kuang,LM,SieberKrauskopf,SC,SS,VTK,
WuWang,ZW}.
It is now well understood that retarded functional differential equations (RFDEs), a class which contains delay-differential equations, behave for the most part like ordinary differential equations on appropriate infinite-dimensional function spaces.  As such, many of the techniques and theoretical results of finite-dimensional dynamical systems have counterparts in the theory of RFDEs.  In particular, versions of the stable/unstable and center manifold theorems in neighborhoods of an equilibrium point exist for RFDEs \cite{HVL}.  Also, techniques for simplifying vector fields via center manifold and normal form reductions have been adapted to the study of bifurcations in RFDEs \cite{FM1,FM2}.

One of the challenges of applying these finite-dimensional techniques to RFDEs lies in the so-called {\it realizability problem}.  This problem stems from the fact that 
the procedure to reduce an RFDE to a center manifold often leads to algebraic restrictions on the nonlinear terms in the center manifold equations. 
Specifically, suppose $B$ is an arbitrary $m\times m$ matrix.  For the
sake of simplicity, suppose additionally that all eigenvalues of $B$
are simple.  
Let $C([-r,0],\mathbb{R})$ be the space of continuous functions
from the interval $[-r,0]$ into $\mathbb{R}$, and for any continuous
function $z$, define $z_t\in
C([-r,0],\mathbb{R})$ as $z_t(\theta)=z(t+\theta)$, $-r\leq\theta\leq 0$.
It is then possible \cite{FM3} to construct a bounded linear operator 
$\mathcal{L}:C([-r,0],\mathbb{R})\longrightarrow\mathbb{R}$ such that the
infinitesimal generator $A$ for the flow associated with the functional
differential equation 
\begin{equation}
\dot{z}(t)=\mathcal{L}\,z_t
\label{linfde1}
\end{equation}
has a spectrum which
contains the eigenvalues of $B$ as a subset. 
Thus, there exists an $m$-dimensional subspace $P$ of
$C([-r,0],\mathbb{R})$ which is invariant for the flow generated by
$A$, and the flow on $P$ is given by the linear ordinary
differential equation (ODE)
\[
\dot{x}=Bx.
\]

Now, suppose (\ref{linfde1}) is modified by the addition of a
nonlinear delayed term
\begin{equation}
\dot{z}(t)=\mathcal{L}\,z_t+az(t-\tau)^2,
\label{nonlinfde1}
\end{equation}
where $a\in\mathbb{R}$ is some coefficient and $\tau\in [0,r]$ is the
delay time.  Then the center manifold theorem for RFDEs \cite{HVL} can be used to
show that the flow for (\ref{nonlinfde1}) admits an $m$-dimensional
locally invariant center manifold on which the dynamics associated
with (\ref{nonlinfde1}) are given by a vector field which, to
quadratic order, is of the form
\begin{equation}
\dot{x}=Bx+ag(x),
\label{realizeode1}
\end{equation}
where $g:\mathbb{R}^m\longrightarrow\mathbb{R}^m$ is a fixed homogeneous
quadratic polynomial which is completely determined by $\mathcal{L}$ and $\tau$,
and $a$ is the same coefficient which appears in (\ref{nonlinfde1}).
We immediately notice that for fixed $\mathcal{L}$ and $\tau$,
(\ref{realizeode1}) has at most one degree of freedom in the quadratic
term, corresponding to the one
degree of freedom in the quadratic term in (\ref{nonlinfde1}).
However, whereas one degree of freedom is sufficient to
describe the general scalar quadratic term involving one delay in
(\ref{nonlinfde1}), it is largely insufficient (if $m>1$) to describe
the general homogeneous quadratic polynomial
$f:\mathbb{R}^m\longrightarrow\mathbb{R}^m$.
Therefore, there exist $m$-dimensional vector fields $\dot{x}=Bx+f(x)$
(where $f$ is homogeneous quadratic) which can not be realized by
center manifold reduction (\ref{realizeode1}) of any RFDE of the form (\ref{nonlinfde1}).  The realizability problem has received considerable attention in the literature \cite{BuonoBelair,ChoiLeBlanc1,ChoiLeBlanc2,FM3,FM4}.

In this paper, we will be interested in a realizability problem for a class of second-order scalar delay-differential equations of the form
\begin{equation}
\ddot{x}(t)+b\dot{x}(t)+ax(t)-F(x(t),\dot{x}(t))=\alpha x(t-\tau)+\beta\dot{x}(t-\tau)+G(x(t-\tau),\dot{x}(t-\tau)),
\label{premaineq}
\end{equation}
where $a$, $b$, $\alpha$ and $\beta$ are real parameters, $\tau>0$ is a delay term, and the nonlinear functions $F$ and $G$ are smooth and vanish at the origin, along with their first order partial derivatives.  This class contains many interesting applications which have been studied in the literature, including Van der Pol's oscillator with delayed feedback \cite{Atay,deOlivera,JiangYuan,WeiJiang,WuWang}, as well as models for stabilization of an inverted pendulum via delayed feedback \cite{SieberKrauskopf}.  

Both the Van der Pol oscillator \cite{WuWang} and the inverted pendulum system \cite{SieberKrauskopf} have been shown to possess points in parameter space where a bifurcation via a non-semisimple triple-zero eigenvalue occurs.  In \cite{SieberKrauskopf}, this bifurcation is in fact characterized as the {\it organizing center} for their model, since it includes in its unfolding Bogdanov-Takens and steady-state/Hopf mode interactions and the associated complex dynamics of these codimension two singularities.  
As far as we are aware, a complete theoretical analysis and classification of all possible dynamics near the non-semisimple triple-zero bifurcation has yet to be done, although a rather thorough investigation was undertaken in \cite{DumortierIbanez}.  Numerical tools are used in \cite{SieberKrauskopf} to illustrate the complexity of this singularity in their model, including many global bifurcations.  It is stated in \cite{SieberKrauskopf} that because of the presence of invariant tori, a full versal unfolding of the triple-zero singularity must include terms other than those appearing at cubic order in their model, and conclude by wondering whether full realizability of the nonlinear normal form for the triple-zero bifurcation is possible for their delay-differential equation.  

Other relevant work includes \cite{CampbellYuan}, where the authors study a class of coupled first-order delay-differential equations which includes (\ref{premaineq}) as a special case (if one writes (\ref{premaineq}) as a first order system), and compute quadratic and cubic normal form coefficients in term of DDE coefficients for both non-semisimple double-zero and triple-zero bifurcations.  Higer-order terms for these normal form are not considered.

In this paper, we will first show that the non-semisimple triple-zero singularity occurs generically in (\ref{premaineq}), and then prove that the full nonlinear normal form for the non-semisimple triple-zero bifurcation, at any prescribed order, can be realized by center manifold normal form reduction of (\ref{premaineq}) for appropriate choices of nonlinear functions $F$ and $G$.  In section 2, we present the functional analytic framework in which we will study this problem.  Section 3 gives a brief summary of the center manifold and normal form procedure for RFDEs which was developed by Faria and Magalh$\tilde{\mbox{\rm a}}$es \cite{FM1,FM2}.  Our main result is stated and proved in section 4.  We end with some concluding remarks in section 5.

\Section{Functional analytic setup}
As mentioned in the introduction, we consider a general class of second order nonlinear differential equations for the real-valued function $x(t)$ of the form 
(\ref{premaineq}), which we rewrite as a first order system
\begin{equation}
\begin{array}{rcl}
\dot{x}(t)&=&y(t)\\
&&\\
\dot{y}(t)&=&-a\,x(t)+\alpha\,x(t-\tau)-b\,y(t)+\beta\,y(t-\tau)+F(x(t),y(t))+G(x(t-\tau),y(t-\tau)).
\end{array}
\label{maineq}
\end{equation}
where $a$, $b$, $\alpha$ and $\beta$ are real parameters, $\tau>0$ is a delay term, and the nonlinear functions $F$ and $G$ are smooth and vanish at the origin, along with their first order partial derivatives.
In many applications, we have $a>0$, so we will assume this condition throughout (although other cases of $a$ can be treated in a similar manner).

The characteristic equation corresponding to (\ref{maineq}) is $P(\lambda)=0$, where
\begin{equation}
P(\lambda)= \lambda^2+b\lambda+a-(\alpha+\beta\lambda)e^{-\lambda\tau}.
\label{char}
\end{equation}
A straightforward computation shows that when
\begin{equation}
\begin{array}{c}
\alpha=a\\
\\
\tau=\tau_0={\displaystyle\frac{\beta}{a}+\frac{\sqrt{\beta^2+2a}}{a}}\\
\\
b=\beta-a\tau_0\\
\\
3\beta^2\neq 2a\,\,\,(\mbox{\rm or equivalently}\,\,\, a\tau_0\neq 3\beta)
\end{array}
\label{param}
\end{equation}
then $P(0)=P'(0)=P''(0)=0$, $P'''(0)\neq 0$, and $P$ has no other roots on the imaginary axis.  Therefore $0$ is a triple eigenvalue for the linearization of (\ref{maineq}) at the origin, with geometric multiplicity one.

For the parameter values (\ref{param}), we write (\ref{maineq}) as 
\begin{equation}
\begin{array}{rcl}
\dot{x}(t)&=&y(t)\\
&&\\
\dot{y}(t)&=&-a\,(x(t)-x(t-\tau_0))+a\tau_o\,y(t)-\beta\,(y(t)-y(t-\tau_0))\\&&\\
&&+F(x(t),y(t))+G(x(t-\tau_0),y(t-\tau_0)).
\end{array}
\label{maineqsys}
\end{equation}
 
Let $C=C\left(  \left[  -\tau_0,0\right]  ,\mathbb{R}^2\right)$ be the
Banach space of continuous functions from $\left[  -\tau_0,0\right]  $ into
$\mathbb{R}^2$ with supremum norm.  We define $z_{t}\in C$ as 
\[z_{t}\left(  \theta\right)
=z\left(  t+\theta\right) =\left(\begin{array}{c} x(t+\theta)\\y(t+\theta)\end{array}\right),-\tau_0\leq\theta\leq 0.
\]
We view (\ref{maineqsys}) as a retarded functional differential equation of the form
\begin{equation}
\dot{z}\left(  t\right)  ={\cal L}\,z_{t}+{\cal F}\left(  z_{t}\right)  ,
\label{y1}
\end{equation}
\noindent where ${\cal L}:C\rightarrow\mathbb{R}^2$ is
the
bounded linear
operator 
\[
{\cal L}\,\phi=\int_{-\tau_0}^{0}\left[  d\eta\left(\theta\right)  \right]  \phi\left(
\theta\right)=\left( \begin{array}{cc}0&1\\-a&a\tau_0-\beta\end{array}\right)\phi(0)+\left(\begin{array}{cc}0&0\\ a&\beta\end{array}\right)\phi(-\tau_0)
\]
and ${\cal F}$ is
the smooth nonlinear function from $C$ into $\mathbb{R}^2$
\[
{\cal F}(\phi)=\left(\begin{array}{c}0\\F(\phi(0))+G(\phi(-\tau_0))\end{array}\right).
\]

Let $A$ be the infinitesimal generator of the flow for the linear
system $\dot{z}={\cal L}\,z_t$, with spectrum $\sigma(A)\supset\,\{0\}$, and
$P$ be the three-dimensional invariant subspace for $A$ associated with the
eigenvalue $0$.  Then it follows that the columns of the matrix
\[
\Phi=\left(\begin{array}{ccc}1&\theta&\frac{1}{2}\theta^2\\0&1&\theta\end{array}\right)
\]
form a basis for $P$.  

In a similar manner, we can define an
invariant space, $P^{\ast},$ to be the generalized eigenspace of the
transposed system, $A^{T}$\ associated with the triple nilpotency having as
basis the rows of the matrix
$\Psi=$\textrm{col}$\left(  \psi_{1},\ldots,\psi_{m}\right)$. Note that the
transposed system, $A^{T}$ is defined over a dual space $C^{\ast}=C\left(
\left[  0,\tau_0\right]  ,\mathbb{R}^2\right),$ and each element of $\Psi$
is included in $C^{\ast}.$ The bilinear form between $C^{\ast}$ and $C$ is
defined as
\begin{equation}
\left(  \psi,\phi\right)  =\psi\left(  0\right)  \phi\left(  0\right)
-\int\limits_{-r}^{0}\int\limits_{0}^{\theta}\psi\left(  \zeta-\theta\right)
\text{ }\left[  d\eta\left(  \theta\right)  \right]  \text{ }\phi\left(
\zeta\right)  \text{ }d\zeta. \label{y11}%
\end{equation}
Note that $\Phi$ and $\Psi$ satisfy $\dot{\Phi}=B\Phi,$
$\dot{\Psi}=-\Psi B,$ where $B$ is the $3\times 3$ matrix
\begin{equation}
B=\left(\begin{array}{ccc}0&1&0\\0&0&1\\0&0&0\end{array}\right).
\label{BBdef}
\end{equation}

We can normalize
$\Psi$ such that $\left(  \Psi
,\Phi\right)  =I$, and we can decompose the space $C$ using the splitting
$C=P\oplus Q$, where the complementary space $Q$ is also invariant for $A$.

Faria and Magalh$\tilde{\mbox{\rm a}}$es \cite{FM1,FM2} show that (\ref{y1}) can be written
as an infinite dimensional ordinary differential equation on the
Banach space 
$BC$ of functions from $[-\tau_0,0]$ into ${\mathbb R}^2$ 
which are uniformly continuous on $[-\tau_0,0)$ and with a jump discontinuity 
at $0$, using a procedure that we will now outline.
Define $X_0$ to be the function
\[
X_0(\theta)=\left\{\begin{array}{lc}
\left(\begin{array}{cc}1&0\\0&1\end{array}\right)&\theta=0\\[0.15in]
\left(\begin{array}{cc}0&0\\0&0\end{array}\right)&-\tau_0\leq\theta<0,
\end{array}\right.
\]
then the elements of $BC$ can be written as $\xi=\varphi+X_0\lambda$, with
$\varphi\in C$ and $\lambda\in {\mathbb R}^2$, so that
$BC$ is identified with $C\times {\mathbb R}^2$.  

Let $\pi:BC\longrightarrow P$ denote the projection
\[
\pi(\varphi+X_0\lambda)=\Phi [(\Psi,\varphi)+\Psi(0)\lambda],
\]
where $\varphi\in C$ and $\lambda\in {\mathbb R}^2$.  
We now decompose $z_t$ in (\ref{y1}) according to the splitting
\[
BC=P\oplus\mbox{\rm ker}\,\pi,
\]
with the property that $Q\subsetneq\,\mbox{\rm ker}\,\pi$, 
and get the following infinite-dimensional ODE system which is equivalent to (\ref{y1}):
\begin{equation}
\begin{array}{rcl}
\dot{u}&=&Bu+\Psi(0)\,{\cal F}(\Phi\,u+v)\\[0.15in]
{\displaystyle\frac{d}{dt}\,v}&=&A_{Q^1}v+(I-\pi)X_0\,{\cal F}(\Phi\,u+v),
\end{array}
\label{projfdep}
\end{equation}
where $u\in {\mathbb R}^3$,
$v\in Q^1\equiv Q\cap C^1$,
($C^1$ is the subset of $C$ consisting
of continuously differentiable functions),
and
$A_{Q^1}$ is the operator from 
$Q^1$ into 
$\mbox{\rm ker}\,\pi$ defined by
\[
A_{Q^1}\varphi=\dot{\varphi}+X_0\,[{\cal L}\,\varphi-\dot{\varphi}(0)].
\]

\Section{Faria and Magalh$\tilde{\mbox{\bf a}}$es normal form}

Consider the formal Taylor expansion of the nonlinear terms $\mathcal{F}$ in (\ref{y1})
\[
{\cal F}(\phi)=\sum_{j\geq 2}\,{\cal F}_j(\phi),\,\,\,\,\,\phi\in\,C,
\]
where ${\cal F}_j(\phi)=V_j(\phi,\ldots,\phi)$, with $V_j$ belonging to the space of
continuous multilinear symmetric maps from
$C\times\cdots\times C$
($j$ times) to $\mathbb{R}^2$. 
If we denote $f_j=(f_j^1,f_j^2)$, where
\[
\begin{array}{rcl}
f_{j}^1(u,v)&=&\Psi(0)\,{\cal F}_j(\Phi\,u+v)\\[0.15in]
f_j^2(u,v)&=&(I-\pi)\,X_0\,{\cal F}_j(\Phi\,u+v),
\end{array}
\]
then (\ref{projfdep}) can be written as
\begin{equation}
\begin{array}{rcl}
\dot{u}&=&{\displaystyle Bu+\sum_{j\geq
    2}\,f_j^1(u,v)}\\[0.15in]
{\displaystyle\frac{d}{dt}\,v}&=&{\displaystyle A_{Q^1}v+\sum_{j\geq
    2}\,f_j^2(u,v)}.
\end{array}
\label{y4}
\end{equation}

It is easy to see that the non-resonance
condition of Faria and Magalh$\tilde{\mbox{\rm a}}$es (Definition (2.15) of \cite{FM2}) holds.  Consequently, using
successively at each order $j$ a near identity change of variables of the form
\begin{equation}
(u,v)=(\hat{u},\hat{v})+U_j(\hat{u})\equiv (\hat{u},\hat{v})+
(U^1_j(\hat{u}),U^2_j(\hat{u})),
\label{nfcv}
\end{equation}
(where $U^{1,2}_j$ are homogeneous degree $j$ polynomials in
the indicated variable, with coefficients respectively in
$\mathbb{R}^3$ and $Q^1$)
system (\ref{y4}) can be put into formal normal form
\begin{equation}
\begin{array}{rcl}
\dot{u}&=&{\displaystyle Bu+\sum_{j\geq
    2}\,g_j^1(u,v)}\\[0.15in]
{\displaystyle\frac{d}{dt}\,v}&=&{\displaystyle A_{Q^1}v+\sum_{j\geq
    2}\,g_j^2(u,v)}
\end{array}
\label{y5}
\end{equation}
such that the center manifold is locally given by $v=0$ and the local
flow of (\ref{y1}) on this center manifold is given by
\begin{equation}
\dot{u}=Bu+\sum_{j\geq 2}\,g_j^1(u,0).
\label{y6}
\end{equation}
The nonlinear terms in (\ref{y6}) are in normal form in the
classical sense with respect to the matrix $B$.

\Section{Realizability of the normal form for the triple-zero nilpotency}

It was shown in \cite{DumortierIbanez,YuYuan} that the classical normal form for the general nonlinear vector field
\[
\left(\begin{array}{c}\dot{u}_1\\\dot{u}_2\\\dot{u}_3\end{array}\right)=\left(\begin{array}{ccc}0&1&0\\0&0&1\\0&0&0\end{array}\right)\,\left(\begin{array}{c}u_1\\u_2\\u_3\end{array}\right)+\left(\begin{array}{c}r_1(u_1,u_2,u_3)\\r_2(u_1,u_2,u_3)\\r_3(u_1,u_2,u_3)\end{array}\right)
\]
(where ${\displaystyle r_j(0,0,0)=0, \frac{\partial r_j}{\partial u_k}(0,0,0)=0, \,\,\mbox{\rm for}\,j,k=1,2,3}$)
is
\[
\begin{array}{ccl}
\dot{u}_1&=&u_2\\
\dot{u}_2&=&u_3\\
\dot{u}_3&=&{\displaystyle \sum_{j\geq 2}\,\left(\,\sum_{i=0}^j\,a_{(j-i),i}\,u_1^{j-i}u_2^i\,+\,u_1^Iu_3\sum_{i=0}^{J}\,b_{N(J-i),i}\,u_1^{J-i}u_3^i\,\right)},
\end{array}
\]
where
\[
\begin{array}{ll}
N=1,\,\,\,J=\frac{1}{2}(j-1),\,\,I=J,\,\,\,\,\,&\mbox{\rm when $j$ is odd},\\&\\
N=2,\,\,\,J=\frac{j}{2}-1,\,\,I=J+1,&\mbox{\rm when $j$ is even}.
\end{array}
\]

Thus, if $B$ is the matrix (\ref{BBdef}), 
and $H^j$ is the space of homogeneous polynomial mappings of degree $j\geq 2$ from $\mathbb{R}^3$ into $\mathbb{R}^3$,
then the homological operator $L_B\equiv Dh(u)\cdot Bu-B\cdot h(u)$ acting on $H^j$ is such that
\begin{equation}
H^j=L_B(H^j)\oplus W_j,
\label{homsplit}
\end{equation}
where $W_j\subset H^j$ is the subspace of dimension
\[
\begin{array}{cl}
{\displaystyle\frac{3j+2}{2}}&\,\,\,\,\mbox{\rm if $j$ is even}\\&\\
{\displaystyle\frac{3j+3}{2}}&\mbox{\rm if $j$ is odd}
\end{array}
\]
spanned by
\[
\left\{\left(\begin{array}{c}0\\0\\u_1^{j-i}u_2^i\end{array}\right)\,,\,i=0,\ldots,j\,\right\}\,\bigcup\,
\left\{\left(\begin{array}{c}0\\0\\u_1^{j-1-i}u_3^{i+1}\end{array}\right)\,,\,
\begin{array}{ll}
i=0,\ldots,\frac{1}{2}(j-2)\,\,\,\,&\mbox{\rm $j$ even}\\&\\
i=0,\ldots,\frac{1}{2}(j-1)&\mbox{\rm $j$ odd}.
\end{array}\,\right\}
\]

Now, if $F(z_1,z_2)$ and $G(z_1,z_2)$  are real-valued functions such that $F$, $G$ and their first-order partial derivatives vanish at the origin, then we may write the Taylor series
\[
F(z_1,z_2)=\sum_{j\geq 2}\,\hat{F}_j(z_1,z_2),\,\,\,\,\,\,\,\,\,\,\,\,\,\,\,\,
G(z_1,z_2)=\sum_{j\geq 2}\,\hat{G}_j(z_1,z_2),
\]
where the $\hat{F}_j$ and $\hat{G}_j$ are homogeneous degree $j$ polynomials.  The first equation in (\ref{y4}) then reduces to
\begin{equation}
\dot{u}=Bu+\sum_{j\geq 2}\left(\begin{array}{c}0\\
\\
\kappa_1\,(\hat{F}_j((u_1,u_2)+v(0))+\hat{G}_j((u_1-\tau_0u_2+\frac{1}{2}\tau_0^2u_3,u_2-\tau_0u_3)+v(-\tau_0)))\\
\\
\kappa_2\,(\hat{F}_j((u_1,u_2)+v(0))+\hat{G}_j((u_1-\tau_0u_2+\frac{1}{2}\tau_0^2u_3,u_2-\tau_0u_3)+v(-\tau_0)))\end{array}\right),
\label{prenf}
\end{equation}
where
\[
\kappa_1=\frac{3(a\tau_0-4\beta)}{2\tau_0(a\tau_0-3\beta)^2},\,\,\,\,\,\,\,\,\,\,\,\,\,\,\,
\kappa_2=\frac{6}{\tau_0^2(a\tau_0-3\beta)}\neq 0.
\]
We note that when $v=0$ in (\ref{prenf}), then 
\[
\begin{array}{ccl}
\hat{F}_j(u_1,u_2)+\hat{G}_j(u_1-\tau_0u_2+\frac{1}{2}\tau_0^2u_3,u_2-\tau_0u_3)&=&
{\cal A}_j(u_1,u_2)+u_1u_3{\cal B}_{j-2}(u_1,u_3)+\\
&&\\
&&u_2u_3{\cal C}_{j-2}(u_1,u_2,u_3)+u_3^2{\cal D}_{j-2}(u_3),
\end{array}
\]
where
\begin{alignat}{2}\label{Adef}
{\cal A}_j(u_1,u_2)=&\,\hat{F}_j(u_1,u_2)+\hat{G}_j(u_1-\tau_0u_2,u_2)\\\label{Bdef}
u_1u_3{\cal B}_{j-2}(u_1,u_3)=&\,\hat{G}_j(u_1+\frac{1}{2}\tau_0^2u_3,-\tau_0u_3)-\hat{G}_j(u_1,0)-\hat{G}_j(\frac{1}{2}\tau_0^2u_3,-\tau_0u_3)\\\notag
u_2u_3{\cal C}_{j-2}(u_1,u_2,u_3)=&\,\hat{G}_j(u_1-\tau_0u_2+\frac{1}{2}\tau_0^2u_3,u_2-\tau_0u_3)-\hat{G}_j(u_1-\tau_0u_2,u_2)-\\\notag
&\,\hat{G}_j(u_1+\frac{1}{2}\tau_0^2u_3,-\tau_0u_3)+\hat{G}_j(u_1,0)\\\notag
u_3^2{\cal D}_{j-2}(u_3)=&\,\hat{G}_j(\frac{1}{2}\tau_0^2u_3,-\tau_0u_3).
\end{alignat}
\begin{lemma}
For a given integer $j\geq 2$, let $\zeta(u_1,u_3)$ be a homogeneous degree $j$ polynomial such that $\zeta(0,u_3)=\zeta(u_1,0)=0$.  Then
there exists a homogeneous polynomial of degree $j$, $\xi(u_1,u_3)$ such that
\begin{equation}
\zeta(u_1,u_3)=\xi(u_1+\frac{1}{2}\tau_0^2u_3,-\tau_0u_3)-\xi(u_1,0)-\xi(\frac{1}{2}\tau_0^2u_3,-\tau_0u_3).
\label{magic}
\end{equation}
\end{lemma}
\noindent
{\bf Proof of lemma:}  If we write
\[
\xi(u_1,u_3)=\sum_{i=0}^j\,\gamma_{j-i,i}\,u_1^{j-i}u_3^i
\]
then a lengthy but straightforward computation shows that
\[
\begin{array}{l}
{\displaystyle \xi(u_1+\frac{1}{2}\tau_0^2u_3,-\tau_0u_3)-\xi(u_1,0)-\xi(\frac{1}{2}\tau_0^2u_3,-\tau_0u_3)=}\\
\\
{\displaystyle\gamma_{j,0}\left(\sum_{k=1}^{j-1}\left(\begin{array}{c}j\\k\end{array}\right)\left(\frac{1}{2}\tau_0^2\right)^k\,u_1^{j-k}u_3^k\right)+
\sum_{i=1}^{j-1}\,\gamma_{j-i,i}\left(
\sum_{k=0}^{j-i-1}\left(\begin{array}{c}j-i\\k\end{array}\right)\left(\frac{1}{2}\tau_0^2\right)^k(-\tau_0)^i\,u_1^{j-i-k}u_3^{i+k}\right)},
\end{array}
\]
where
\[
\left(\begin{array}{c}j\\k\end{array}\right)=\frac{j!}{k!(j-k)!}.
\]
Now, since $\zeta(0,u_3)=\zeta(u_1,0)=0$, we have that
\[
\zeta(u_1,u_3)=\sum_{i=1}^{j-1}\,\epsilon_{j-i,i}u_1^{j-i}u_3^i.
\]
Thus, we see for example that we may solve (\ref{magic}) by arbitrarily setting $\gamma_{j,0}=0$, $\gamma_{0,j}=0$, and choosing $\gamma_{j-i,i}$, $i=1,\ldots,j-1$ such that the following triangular linear algebraic system is satisfied
\[
\begin{array}{c}
{\displaystyle\left(\begin{array}{c}j-1\\0\end{array}\right)(-\tau_0)\,\gamma_{j-1,1} = \epsilon_{j-1,1}}\\
\\
{\displaystyle\left(\begin{array}{c}j-1\\1\end{array}\right)\left(\frac{1}{2}\tau_0^2\right)(-\tau_0)\,\gamma_{j-1,1}+
\left(\begin{array}{c}j-2\\0\end{array}\right)(-\tau_0)^2\,\gamma_{j-2,2}=\epsilon_{j-2,2}}\\
\\
\vdots\\
\\
{\displaystyle\left(\begin{array}{c}j-1\\j-2\end{array}\right)\left(\frac{1}{2}\tau_0^2\right)^{j-2}(-\tau_0)\,\gamma_{j-1,1}+\hdots+
\left(\begin{array}{c}1\\0\end{array}\right)(-\tau_0)^{j-1}\,\gamma_{1,j-1}=\epsilon_{1,j-1}}
\end{array}
\]
This ends the proof of the lemma.
\hfill\qed

Now, recalling the splitting (\ref{homsplit}), let $\Theta_j(u_1,u_2,u_3)$ be a homogeneous degree $j$ polynomial such that $\Theta_j\in W_j$.
We may write
\[
\Theta_j(u)=\left(\begin{array}{c}0\\0\\q_j(u_1,u_2)\end{array}\right)+\left(\begin{array}{c}0\\0\\u_1u_3s_{j-2}(u_1,u_3)\end{array}\right),
\]
where $q_j$ is a homogeneous polynomial of degree $j$, and $s_{j-2}$ is a homogeneous degree $j-2$ polynomial.

The degree $j$ term in (\ref{prenf}) for $v=0$ can be written as
\[
\begin{array}{l}
\left(\begin{array}{c}0\\
\\
\kappa_1\,(\hat{F}_j(u_1,u_2)+\hat{G}_j(u_1-\tau_0u_2+\frac{1}{2}\tau_0^2u_3,u_2-\tau_0u_3))\\
\\
\kappa_2\,(\hat{F}_j(u_1,u_2)+\hat{G}_j(u_1-\tau_0u_2+\frac{1}{2}\tau_0^2u_3,u_2-\tau_0u_3))\end{array}\right)=\\
\\
R_j(u_1,u_2,u_3)+\left(\begin{array}{c}0\\0\\\kappa_2\,{\cal A}_j(u_1,u_2)\end{array}\right)+
\left(\begin{array}{c}0\\0\\\kappa_2\,u_1u_3{\cal B}_{j-2}(u_1,u_3)\end{array}\right)\end{array}
\]
where $R_j$ is in the range of the homological operator, $R_j\subset L_B(H^j)$, and ${\cal A}_j$ and ${\cal B}_{j-2}$ are as in (\ref{Adef}) and (\ref{Bdef}).
Using the previous lemma, we know that if we choose $\hat{G}_j$ such that $\kappa_2\,{\cal B}_{j-2}=s_{j-2}$, and then set $\kappa_2\,\hat{F}_j(u_1,u_2)=q_j(u_1,u_2)-\kappa_2\,\hat{G}_j(u_1-\tau_0u_2,u_2)$, then
\begin{equation}
\Theta_j(u)=\left(\begin{array}{c}0\\0\\\kappa_2\,{\cal A}_j(u_1,u_2)\end{array}\right)+
\left(\begin{array}{c}0\\0\\\kappa_2\,u_1u_3{\cal B}_{j-2}(u_1,u_3)\end{array}\right).
\label{theteq}
\end{equation}

We can now state and prove the following realizability theorem:
\begin{thm}
Given an integer $\ell\geq 2$ and a polynomial vector field on $\mathbb{R}^3$ of the form
\begin{equation}
\dot{u}=Bu+\sum_{j=2}^{\ell}\,w_j(u),
\label{nf2}
\end{equation}
where $B$ is the matrix (\ref{BBdef}) and $w_j\in W_j$ as in (\ref{homsplit}), there exist polynomial functions $F$ and $G$ in (\ref{maineqsys}) such that the Faria and 
Magalh$\tilde{\mbox{\it a}}$es center manifold and normal form reduction (\ref{y6}) of (\ref{prenf}) up to order $\ell$ is (\ref{nf2}).
\label{mainthm}
\end{thm}
\noindent
{\bf Proof:}  Applying successively at each order $j$ (from $j=2$ to $j=\ell$) a near identity change of variables of the form (\ref{nfcv}), and setting $v=0$, we transform (\ref{prenf}) into
\begin{equation}
\dot{u}=Bu+\sum_{j=2}^{\ell}\,\left[\left(\begin{array}{c}0\\0\\\kappa_2\,{\cal A}_j(u_1,u_2)\end{array}\right)
+\left(\begin{array}{c}0\\0\\\kappa_2\,u_1u_3{\cal B}_{j-2}(u_1,u_3)\end{array}\right)+\Lambda_j(u_1,u_2,u_3)\right]+O(|u|^{\ell+1}),
\label{postnf}
\end{equation}
where $\Lambda_2(u_1,u_2,u_3)=0$ and for $j\geq 3$, $\Lambda_j\in W_j$ is an extra contribution to the terms of order $j$ coming from the transformation of the lower order $(<j)$ terms.
Therefore, we set
\[
\Theta_j(u)=w_j(u)-\Lambda_j(u),\,\,\,\,\,j=2,\ldots,\ell,
\]
and use $(\ref{theteq})$ to conclude that  the truncation of (\ref{postnf}) at order $\ell$ is (\ref{nf2}).
\hfill\qed

\Section{Conclusion}

In this paper, we have solved the realizability problem for the normal form of the non-semisimple triple-zero singularity in a class of delay differential equations (\ref{premaineq}) which includes delayed Van der Pol oscillators, as well as certain models for the control of an inverted pendulum as special cases.    It is apparent from the complexity of the dynamics of (\ref{premaineq}) near the triple-zero nilpotency reported in previous work \cite{CampbellYuan,DumortierIbanez,SieberKrauskopf} that high-order normal forms will be required for a complete classification of this singularity.  Although such a complete classification of the dynamics near a triple-zero nilpotency is beyond the scope of this paper, our results allow us to conclude that the full range of complexity of this singularity is attainable within the class of delay-differential equations (\ref{premaineq}).

Although we have not done so, we believe that the results in this paper could be suitably generalized to studying realizability of higher order nilpotencies in higher-order scalar delay-differential equations such as
\[
x^{(n)}(t)+\sum_{j=0}^{n-1}\,a_j\,x^{(j)}(t)-F(x(t),\ldots,x^{(n-1)}(t))=\sum_{j=0}^{n-1}\,\alpha_j\,x^{(j)}(t-\tau)+G(x(t-\tau),\ldots,x^{(n-1)}(t-\tau))
\]
where $n\geq 3$ is an integer.

\vspace*{0.25in}
\noindent
{\Large\bf Acknowledgments}

\vspace*{0.2in}
This research is partly supported by the
Natural Sciences and Engineering Research Council of Canada in the
form of a Discovery Grant.


\begin{thebibliography}{10}

\bibitem{Atay}
F.M.~Atay.
\newblock Van der Pol's oscillator under delayed feedback.
\newblock {\em J. Sound and Vibration} {\bf 218}, (1998) 333--339.

\bibitem{BBL} 
A.~Beuter, J.~B\'elair and C.~Labrie. 
\newblock Feedback and delays in neurological diseases : a modeling study 
using dynamical systems.
\newblock {\em Bulletin Math. Biology} {\bf 55}, (1993) 525--541.

\bibitem{BuonoBelair}
P.-L.~Buono and J.~B\'elair.
\newblock Restrictions and unfolding of double Hopf bifurcation in functional differential equations.
\newblock {\em  J. Differential Equations} {\bf 189}, (2003) 234--266.

\bibitem{CampbellYuan}
S.A.~Campbell and Y.~Yuan.
\newblock Zero singularities of codimension two and three in delay differential equations.
\newblock {\em Nonlinearity} {\bf 21}, (2008) 2671--2691.

\bibitem{ChoiLeBlanc1}
Y.~Choi and V.G.~LeBlanc.
\newblock Toroidal normal forms for bifurcations in retarded functional differential equations I: multiple Hopf and transcritical/multiple Hopf interaction.
\newblock {\em J. Differential Equations} {\bf 227}, (2006) 166--203.

\bibitem{ChoiLeBlanc2}
Y.~Choi and V.G.~LeBlanc.
\newblock Toroidal normal forms for bifurcations in retarded functional differential equations II: saddle-node/multiple Hopf interaction.
\newblock {\em Dynamics of Continuous, Discrete and Impulsive Systems A} {\bf 15} (2008), 251--276.

\bibitem{deOlivera}
J.C.F.~de Oliveira.
\newblock Oscillations in a van der Pol equation with delayed argument.
\newblock {\em J. Math. Anal. Appl.} {\bf 275}, (2002) 789--803.

\bibitem{DumortierIbanez}
F.~Dumortier and S.~Ib\`a$\tilde{\mbox{\rm n}}$ez.
\newblock Nilpotent singularities in generic 4-parameter families of 3-dimensional vector fields.
\newblock {\em J. Differential Equations} {\bf 127}, (1996) 590--647.
 
\bibitem{FM1}
T.~Faria and L.T.~Magalh$\tilde{\mbox{\rm a}}$es.
\newblock Normal forms for retarded functional differential equations and applications to Bogdanov-Takens singularity. 
\newblock {\em J. Differential Equations} {\bf 122}, (1995) 201--224.	

\bibitem{FM2}
T.~Faria and L.T.~Magalh$\tilde{\mbox{\rm a}}$es.
\newblock Normal forms for retarded functional differential equations with parameters and applications to Hopf bifurcation.
\newblock {\em J. Differential Equations} {\bf 122}, (1995) 181--200.

\bibitem{FM3}
T.~Faria and L.T.~Magalh$\tilde{\mbox{\rm a}}$es. 
\newblock Realization of ordinary differential equations by retarded functional differential equations in neighborhoods of equilibrium points.
\newblock {\em Proc. Roy. Soc. Edinburgh Ser. A} {\bf 125}, (1995) 759--776.

\bibitem{FM4}
T.~Faria and L.T.~Magalh$\tilde{\mbox{\rm a}}$es.
\newblock Restrictions on the possible flows of scalar retarded functional differential equations in neighborhoods of singularities.
\newblock {\em J. Dynam. Differential Equations} {\bf 8}, (1996) 35--70.

\bibitem{HVL}
J.K.~Hale and S.M.~Verduyn Lunel.
\newblock Introduction to Functional Differential Equations, 
\newblock Appl. Math. Sci., vol. 99, Springer, New York, 1993.

\bibitem{HFEKGG}
T.~Heil, I.~Fischer, W.~Els\"{a}\mbox{\ss}er, B.~Krauskopf, K.~Green
and A.~Gavrielides.
\newblock Delay dynamics of semiconductor lasers with short external
cavities: Bifurcation scenarios and mechanisms.
\newblock {\em Phys. Rev. E} {\bf 67}, (2003) 066214-1--066214-11.

\bibitem{JiangYuan}
W.~Jiang and Y.~Yuan.
\newblock Bogdanov-Takens singularity in Van der Pol's oscillator with delayed feedback.
\newblock {\em Phys. D} {\bf 227}, (2007) 149--161.

\bibitem{Kuang}
Y.~Kuang.
\newblock {\em Delay differential equations with applications in population 
dynamics.} 
\newblock Mathematics in Science and Engineering, 191. Academic Press,  
Boston, (1993).

\bibitem{LM} 
A.~Longtin and J.G.~Milton. 
\newblock Modeling autonomous oscillations in the human pupil light 
reflex using nonlinear delay-differential equations.
\newblock {\em Bulletin Math. Biology} {\bf 51}, (1989) 605--624.

\bibitem{SieberKrauskopf}
J.~Sieber and B.~Krauskopf.
\newblock Bifurcation analysis of an inverted pendulum with delayed feedback control near a triple-zero eigenvalue singularity.
\newblock {\em Nonlinearity} {\bf 17}, (2004) 85--103.

\bibitem{SC}
E.~Stone and S.A.~Campbell.
\newblock Stability and bifurcation analysis of a nonlinear DDE model
for drilling.
\newblock {\em J. Nonlinear Sci.} {\bf 14}, (2004) 27--57.

\bibitem{SS} M.J. Suarez and P.L. Schopf. 
\newblock A Delayed Action Oscillator for ENSO. 
\newblock
{\em J. Atmos. Sci.}  {\bf 45}, (1988), 3283--3287.

\bibitem{VTK}
A.G.~Vladimirov, D. Turaev and G. Kozyreff.
\newblock Delay differential equations for mode-locked semiconductor
lasers.
\newblock {\em Optics Letters} {\bf 29}, (2004) 1221-1223.

\bibitem{WeiJiang}
J.~Wei and W.~Jiang.
\newblock Stability and bifurcation analysis in Van der Pol's oscillator with delayed feedback.
\newblock {\em J. Sound and Vibration} {\bf 283}, (2005) 801--819.

\bibitem{WuWang}
X.~Wu and L.~Wang.
\newblock Zero-Hopf bifurcation for van der Pol's oscillator with delayed feedback.
\newblock {\em J. Comput. Appl. Math.} {\bf 235}, (2011) 2586--2602.

\bibitem{YuYuan}
P.~Yu and Y.~Yuan.
\newblock The simplest normal forms associated with a triple zero eigenvalue of indices one and two.
\newblock {\em Nonlin. Anal.} {\bf 47}, (2001) 1105--1116.

\bibitem{ZW}
C.~Zhang and J.~Wei.
\newblock Stability and bifurcation analysis in a kind of business
cycle model with delay.
\newblock {\em Chaos Solitons Fractals} {\bf 22}, (2004) 883--896.



\end{thebibliography}
\end{document}